\documentclass[12pt]{article}
\usepackage{amsthm,amstext,amscd,latexsym,amsmath,amsfonts,amssymb}
\begin{document}
\title{Particles, Generalized Statistics and Categories
\thanks{The work is partially sponsored
by Polish Committee for Scientific Research (KBN) under Grant
No 5P03B05620.}}
\author{W\l adys\l aw Marcinek\\
Institute of Theoretical Physics\\
University of Wroc{\l}aw\\
Poland }
\date{}
\maketitle
\begin{abstract}
Interacting systems of particles with generalized statistics are
considered on both classical and quantum level. It is shown that
all possible quantum states and corresponding processes can be
represented in terms of certain specific categories. The
corresponding Fock space representation is discussed. The problem
of existence of well--defined scalar product is considered. It is
shown that commutation relations corresponding to a system with
generalized statistics can be constructed from relations
corresponding to Bolzman statistics.
\end{abstract}

\section{\hspace{-4mm}.\hspace{2mm}Introduction}

Let us assume that we have a system of $n$--identical and spinless particles moving
on a smooth manifold $\mathcal{M}$ without boundary, of dimension $d\geq 2$. All
possible classical trajectories of the particle system are path in the $n$-th
Cartesian product $\mathcal{M}^{\times n}$. It is natural to assume that two or more
particles can not occupy the same position. We also assume that our particles are
identical. This means that the configuration space for the system of particles is
\begin{equation} Q_{n}(\mathcal{M}) = {\left(\mathcal{M}^{\times n}\backslash D_n \right) }/{S_{n}},
\end{equation} where $D$ is the subset of the Cartesian product $\mathcal{M}^{\times n}$ on
which two or more particles occupy the same position and $S_{n}$ is the symmetric
group. Let us denote by $P_m\mathcal{M}$ the space of all homotopy classes of
trajectories which starts at certain point $m_0 \in \mathcal{M}$ and end at arbitrary
point $m \in \mathcal{M}$. The homotopy classes $P_{m_0}\mathcal{M}$ of all paths
which start at $m_0$ and end at the same point $m_0$, (i.e. a loop space) can be
naturally endowed with a group structure. This group is known as the fundamental
group of $\mathcal{M}$ at $m_0$ and is denoted by $\pi_1 (\mathcal{M}, m_0)$. One can
prove that under some assumptions all fundamental groups $\pi_1 (\mathcal{M}, m_0)$
corresponding to different base point $m_0$ are isomorphic one to other. Hence we can
define the fundamental group $\pi_1 (\mathcal{M})$ of $\mathcal{M}$ as $\pi_1
(\mathcal{M}, m_0)$ for arbitrary $m_0$. We denote by $\pi_* :
\pi_1(\tilde{\mathcal{M}}) \rightarrow \pi_1(M)$ the homomorphism induced by $\pi$.
\begin{em}
The fundamental group $\pi_{1}\left( Q_{n}(\mathcal{M})\right)$ of the space
$Q_{n}(\mathcal{M})$ is said to be the $n$--string braid group on $\mathcal{M}$ and is
denoted by $B_{n} (\mathcal{M})$, i.e.
\begin{equation} B_{n} (\mathcal{M}) := \pi_{1}\left( Q_{n}(\mathcal{M})\right). \end{equation}
\end{em}
Note that elements of $B_{n}(\mathcal{M})$ are loops in $Q_n$. It is easy to see that
all loops in $B_{n}(\mathcal{M})$ can be naturally divided into two classes. The
first class corresponds to noncontractible loops with any interchanges of particles.
This class describe the topology of the base space $\mathcal{M}$. The second class
contain loops which corresponds to all mutual interchanges of particles one by other.
It describes the statistics of a given system of particles. Let us consider this class
in more details. Let $\Sigma_i = (i, i+1)\in S_n$ be the transposition which
describes the interchange of the $i$-th particle with the $i+1$-th one. Let us denote
by $\xi_i$ the path in $\mathcal{M}^{\times n} \backslash D_n$ which realize this
interchange. We assume that for the path $\xi_i$ we have the following
parametrization \begin{equation} \xi_i : t \in [0, 1] \rightarrow \xi_i(t) =
[(\xi_i)_1,\ldots,(\xi_i)_n](t).
\end{equation} We assume in addition that all $(\xi_i)_k$ for $k$ not in $\{i, i+1\}$
are fixed, the loop $[(\xi_i)_i, (\xi_i)_{i+1}](t) \subset \mathcal{M}$ is
counterclockwise and all particles $k$ not in $\{i, i+1\}$ are outside of this loop.
The projection of $\xi_i$ on $Q_n$ is a loop in $Q_n$. The homotopy class of this
loop is denoted by $\Sigma_i$. One can see that elements $\{s_i : i = 1,\ldots,n-1\}$
satisfy two following conditions
\begin{equation} \begin{array}{lll} s_i s_{i+1} s_i = s_{i+1} s_i s_{i+1}&
\mbox{for}&j = 1,\ldots,n-2,\\
s_i s_j = s_j s_i &\mbox{for}&|i-j| > 2. \end{array} \end{equation} This means that
$s_i$ for $(i = 1,\ldots,n-1)$ generate a group which is denoted by $\Sigma_n
(\mathcal{M})$.

Let $\tilde{Q}_n$ be the space which cover the configuration space $Q_n$. The above
definition means that we have the relation \begin{equation} \Sigma_n (\mathcal{M}) =
\pi_{\ast}[\pi_1(\tilde{Q} _n (\mathcal{M}))]. \end{equation} Note that the group
$\Sigma_n(\mathcal{M})$ is a subgroup of $B_n(\mathcal{M})$ and is an extension of the
symmetric group $S_n$ describing the interchange process of two arbitrary
indistinguishable particles. It is obvious that the statistics of the given system of
particles is determined by this group, see \cite{Wu,I}. A system of $n$--identical
and spinless particles moving on a manifold $\mathcal{M}$ and equipped with statistics
described by the group $\Sigma_i(\mathcal{M})$ is said to be a particle system with
the braid-group statistics.

We must indicate that the presented above braid group approach to quantum  arbitrary
statistics is in our opinion not adequate for the study of charged particles
interacting with an external quantum field. First of all the statistics is determined
by the interchange of two different particles. In order to do such interchange we
need two particles and a place, a suitable space.

The statistics of a single particle has no sense. Hence we need a new and more
general approach to the concept of statistics. In this paper we are going to study a
system of charged particles moving under influence of a quantum field. A generalized
statistics is considered as a result of interaction of charged particles with quantum
field. Our discussion is based on the algebraic formalism previously considered in
\cite{mad,qweyl,mco,top,castat,quo}. The short summary of this paper has been
published in the Proceedings of the XII-th Max Born Symposium \cite{gqs}. The
fundamental assumption is that every charged particle is transform under interaction
into a system consisting a charge and  quanta of the field. Such system behaves like
free particles moving in certain effective space. The so--called cross statistics is
discussed. This statistics is described by an operator $T$. The existence of exchange
braid statistic is related to the existence of additional operator $B$ and
corresponding consistency conditions. For the description of our system we need a
solution of these conditions. It is shown that these solutions are related to a
construction of a category of states.

\section{\hspace{-4mm}.\hspace{2mm}Interaction as coaction}

Let us consider for example a charged particle moving on a space $\mathcal{M}$ under
influence of singular magnetic field \cite{sin}. We assume that $\mathcal{M}$ is a
path-connected topological space with a base point $m_0$. All possible classical
trajectories of the particle are path in $\mathcal{M}$. We denote by $P_m\mathcal{M}$
the space of all homotopy classes of paths which starts at $m_0$ and end at arbitrary
point $m \in \mathcal{M}$. It is obvious that the union $\cup_m P_m\mathcal{M} =
P\mathcal{M}$ is a covering space $P = (P\mathcal{M}, \pi, \mathcal{M})$. The
projection $\pi : P\mathcal{M} \rightarrow \mathcal{M}$ is given by
\begin{equation}
\begin{array}{ccc} \pi (\xi) = m &\mbox{iff}&\; \xi \in P_m \mathcal{M} . \end{array} \end{equation} The
homotopy class $P_{m_0}\mathcal{M}$ of all paths which start at $m_0$ and end at the
same point $m_0$, (i.e. a loop space) forms the fundamental group $\pi_1
(\mathcal{M}, m_0)$ of $\mathcal{M}$ at $m_0$. Generators of the fundamental group are
denoted by $\Sigma_i$. In our case we have $\mathcal{M} = S^1 \times \ldots \times
S^1$ and the fundamental group is
\begin{equation} \pi_1 (\mathcal{M}, m_0) = Z \oplus...\oplus Z \;\;\; \mbox{(N-sumands)}.
\end{equation} Let $G$ be a subgroup of the fundamental group $\pi_1 (\mathcal{M}, m_0)$. We assume that
quantum states of magnetic field can be represented as linear combinations of
elements (over a field $\mathbb{C}$ of complex numbers) of the group $G$. It is known
such that linear combinations of elements of certain group $G$ form a group algebra
$H := \mathbb{C} G$. It is also known that there is a Hopf algebra structure defined
on the group algebra $H$. Let us denote by $E$ a Hilbert space of quantum states of
particle which is not coupled to magnetic field. Quantum states of particle coupled
to our singular magnetic field is described by the tensor product $E \otimes H$. Every
attaching of magnetic flux to the particle moving in singular magnetic field can be
represented by certain coaction $\rho_E$ of the Hopf algebra $H$ on the space $E$,
i.e. by a linear mapping \begin{equation} \rho_E : E \rightarrow E \otimes H,
\end{equation} which define a (right-) $H$-comodule structure on $E$. Note that if
$E$ is a $H$-comodule, where $H = \mathbb{C} G$, then $E$ is also a $G$-graded vector
space, i.e
\begin{equation} E = \bigoplus\limits_{\alpha \in G} \ E_{\alpha}. \end{equation} Note that
the family of all $H$-comodules for a given Hopf algebra $H$ forms a monoidal
category $\mathcal{C} = \mathcal{M}^H$. In this way we obtain a category of comodules
for the description of particles in singular magnetic field. This means that there is
a closed relation between interacting particle system and categories. Obviously we can
generalize the above particular case and study interactions and quantum processes in
terms of categories in an general way.

\section{\hspace{-4mm}.\hspace{2mm}Categories and interactions}

We are going here to study systems of charged particles with certain dynamical
interaction with a quantum field. It is natural to expect that some new and specific
quantum states of the system have appear as a result of interaction. We would like to
describe all such states. It is natural to assume that there is a specific category
$\mathcal{C}$ which represent all possible quantum states of interacting systems of
particles and corresponding processes. The category $\mathcal{C}$ is said to be a {\it
category of states}.

Our fundamental assumption is that every charged particle is transform under
interaction into a system consisting a charge and $N$--species of quanta of the
field. A system which contains a charge and certain number of quanta as a result of
interaction with the quantum field is said to be a {\it dressed particle}. Next we
assume that every dressed particle is a composite object equipped with an internal
structure. Obviously the structure of dressed particles is determined by the
interaction with the quantum field. We describe a dressed particle as a non-local
system which contains $n$ centers (vertexes). Two systems with $n$ and $m$ centers,
respectively, can be "composed" into one system with $n+m$ centers. All centers as
members of a given system behave like free particles moving on certain effective
space. Every center is also equipped with ability for absorption and emission of
quanta of the intermediate field. A centre dressed with a single quantum of the field
is said to be a {\it quasiparticle}. In our approach a center equipped with two
quanta forms a system of two quasiparticles.

A center with an empty place for a single quantum is said to be a {\it quasihole}. A
centre which contains any quantum is said to be {\it neutral}. A neutral center can
be transform into a quasiparticle or a quasihole by an absorption or emission process
of single quantum, respectively. In this way the process of absorption of quanta of
quantum field by a charged particle is equivalent to a creation of quasiparticles and
emission -- to annihilation of quasiparticles. Note that there is also the process of
mutual annihilation of quasiparticles and quasiholes.

Quasiparticles and quasiholes as components of certain dressed particle have also
their own statistics. It is interesting that there is a statistics of new kind,
namely a {\it cross statistics}. This statistics is determined by an exchange process
of quasiholes and quasiparticles. Note that the exchange is not a real process but an
effect of interaction. Such exchange means annihilation of a quasiparticle on certain
place and simultaneous creation of quasihole on an another place.

Summarizing we assume that for a system of quasiparticles and quasiholes we have the
following possible processes: (i) composition, (ii) creation and annihilation, (iii)
mutual exchanges of quasiholes and quasiparticles, (iv) exchange of quasiparticles or
quasiholes itself. We are going to use the concept of monoidal categories in order to
describe these processes. Let $\mathcal{C} \equiv \mathcal{C}(\otimes, k)$ be a
monoidal category. This means that there are a collection $\mathcal{O}
b(\mathcal{C})$ of objects, a collection $hom(\mathcal{C})$ of arrows (morphisms), an
identity object $k$ (a field) and a monoidal operation $\otimes :{\cal C}\times{\cal
C}\rightarrow {\cal C}$ satisfying some known axioms, see \cite{ML} for details. The
collection $hom(\mathcal{C})$ is the union of mutually disjoint sets $hom(\mathcal{U},
\mathcal{V})$ of arrows $f : \mathcal{U}\rightarrow\mathcal{V}$ from $\mathcal{U}$ to
$\mathcal{V}$ defined for every pair of objects $\mathcal{U}, \mathcal{V} \in
\mathcal{O} b(\mathcal{C})$. It may happen that for a pair $\mathcal{U},
\mathcal{V}\in\mathcal{O} b(\mathcal{C})$ the set $hom(\mathcal{U}, \mathcal{V})$ is
empty. The associative composition of morphisms is also defined. The monoidal
operation $\otimes :{\cal C}\times{\cal C}\rightarrow {\cal C}$ is a bifunctor which
has a two-sided identity object $k$. The category $\mathcal{C}$ can contain some
special objects like algebras, coalgebras, modules or comodules, etc...

Our fundamental assumption is that all possible quantum processes are represented as
arrows of certain monoidal category $\mathcal{C} = \mathcal{C}(\otimes, k)$. In our
case $k\equiv\mathbb{C}$ is the field of complex numbers. If $f : \mathcal{U}
\rightarrow \mathcal{V}$ is an arrow from $\mathcal{U}$ to $\mathcal{V}$, then the
object $\mathcal{U}$ represents physical objects before interactions and $\mathcal{V}$
represents possible results of interactions. We assume that different objects of the
monoidal category $\mathcal{C}$ describe physical objects of different nature, charged
particles, quasiparticles or different species of quanta of an external field, etc...
Let $\mathcal{U}$ be an object of the category $\mathcal{C}$, then the object
$\mathcal{U}^{\ast}$ corresponds for antiparticles, holes or quasiholes or dual
field, respectively. If $\mathcal{U}$ and $\mathcal{V}$ are two different objects of
the category $\mathcal{C}$, then $\mathcal{U}\otimes\mathcal{V}$ is also an object of
the category, it represents a composite quantum system composed from object of
different nature.

Let $\mathcal{U}$ and $\mathcal{V}$ be two different objects of $\mathcal{M}$. If for
example $\mathcal{U}$ represents charged particles and $\mathcal{V}$ -- a quantum
field, then the product $\mathcal{U}\otimes\mathcal{V}$ describes the composite
system containing both particles and the quantum field. Observe that the arrow
$\mathcal{U}\rightarrow\mathcal{U}\otimes\mathcal{V}$ describes the process of
absorption and the arrow $\mathcal{U}\otimes\mathcal{V}\rightarrow\mathcal{U}$
describes the process of emission.

If $\mathcal{B}$ is a coalgebra in $\mathcal{C}$, then the result of comultiplication
$\triangle : \mathcal{B}\rightarrow\mathcal{B}\otimes\mathcal{B}$ represents a
composite object composed from copies of objects of the same nature. in other words
the arrow $\triangle : \mathcal{C}\rightarrow\mathcal{C}\otimes\mathcal{C}$
represents a process of making a copy of an object.

Let $\mathcal{A}$ be an algebra in the category $\mathcal{M}$. The multiplication $m :
\mathcal{A}\otimes\mathcal{A}\rightarrow\mathcal{A}$ is a morphism in this category.
In our interpretation it represents the creation process of a single object from a
composite system of objects of the same species.

Let $\mathcal{A}$ be an algebra and $\mathcal{B}$ be a coalgebra in $\mathcal{C}$. If
we have a bilinear pairing $<-,-> : \mathcal{B}\otimes\mathcal{A}\rightarrow k$ such
that
\begin{equation} <\triangle f , s \otimes t> = <f , s\cdot t>, \end{equation} where
$f\in\mathcal{B}$ and $s, t\in\mathcal{A}$. then we say that $\mathcal{B}$ and
$\mathcal{A}$ are in duality. In our physical interpretation this duality means that
the annihilation process of quasiparticles as components of two independent systems is
equivalent to their annihilation as components of one bigger system!

The cross symmetry
$\Psi_{\mathcal{U}^{\ast},\mathcal{V}}:\mathcal{U}^{\ast}\otimes\mathcal{V}\rightarrow\mathcal{V}\otimes
\mathcal{U}^{\ast}$ corresponds for exchanging process of quasiparticles and
quasiholes. If $\mathcal{V}\equiv\mathcal{U}$, then the pairing
$g_{\mathcal{U}}:\mathcal{U}^{\ast}\otimes\mathcal{U}\rightarrow\mathbb{C}$ describes
the process of annihilation of a pair, quasiparticle and quasihole. In this way we can
characterize in general a category $\mathcal{C}$ representing all possible quantum
processes of certain nature.

Let $\mathcal{U}$ represents physical objects before interactions and $\mathcal{V}$
represents possible results of interactions. If results of a given interaction can be
obtained in a few different ways, then we assume that the only difference is in the
phase. We identify all objects with different phase factor. Such identification is a
natural isomorphism in the category. In this case we say that the process of
interaction is consistent. If every possible process which can be represented as
certain sequence arrows and objects in the category $\mathcal{C}$ is consistent, then
this category is said to be consistent. The construction and classification of
categories for the description of quantum possible processes in a consistent way is a
problem.


\section{\hspace{-4mm}.\hspace{2mm}Quantum states and generalized
statistics}

In this section we are going to formulate conditions which form a starting point for
the construction of category representing particles dressed with some quanta. Let us
consider a single center of interaction which can be transform into system which
contains a charged particle dressed in one quantum in $N$ different way, i. e. a
system of $N$ quasiparticles. In our approach such system is represented by as a
finite set of elements \begin{equation} Q := \{x^i : i = 1,\ldots, N<\infty\}
\end{equation} which form a basis for a finite--dimensional Hilbert space $E$ over a
field of complex numbers $\mathbb{C}$. A center with empty places (quasiholes)
correspond to the basis \begin{equation} Q^{\ast} := \{x^{\ast i} : i = N, N-1,\ldots
, 1\}.
\end{equation} for the complex conjugate space $E^{\ast}$.

The process of "composition" of quasiparticles and quasiholes is described by a
tensor product of spaces $E$ and $E^{\ast}$, respectively. In this way the product
$\underbrace{E\otimes\cdots\otimes E}_n$ represents a system of $n$ centers which can
be equipped with several quanta in $n N$ different ways. The product
$\underbrace{E^{\ast}\otimes\cdots\otimes E^{\ast}}_n$ represents a system of $n$
centers with empty places and an arbitrary mixed product of spaces $E$ and $E^{\ast}$
represents centers equipped with quanta and empty places.

The field $\mathbb{C}$ of complex numbers represent the state without centers of
interaction, i. e. a free and undressing particle in a vacuum state. The pairing $g_E
: E^{\ast} \otimes E \rightarrow \mathbb{C}$ and the corresponding scalar product is
given by
\begin{equation} g_E (x^{\ast i}\otimes x^j) = \langle x^{i}|x^j \rangle := \delta^{ij}. \end{equation} It
represent the annihilation of  empty places by corresponding quanta, i. e.
annihilation of pairs, quasiholes by quasiparticles .

The cross statistics is described by an operator $T$ called an {\it elementary cross}
or {\it twist} \cite{qstat,bma}. This operator is linear, invertible and Hermitian.
It is given by its matrix elements $T : E^{\ast}\otimes E\rightarrow E \otimes
E^{\ast}$
\begin{equation} \begin{array}{c} T(x^{\ast i}\otimes x^j) = \Sigma \ T^{ij}_{kl} x^k \otimes x^{\ast
l}. \label{cross} \end{array} \end{equation} We assume here that both spaces $E$ and
$E^{\ast}$, the $\ast$--operation, the pairing $g_E$ and operators $T$  are the
starting point for a construction of our category $\mathcal{C}$. The most simple
example is provided by the category $\mathcal{C} (E, E^{\ast}, g_{E}, T)$ which
consists of all possible tensor product of spaces $E$ and $E^{\ast}$ and their direct
sums. The $\ast$--operation, the pairing $g_E$ and the cross $T$ can be extended to
the whole category \cite{castat}.

One can also add some quotients and a braid symmetry \cite{castat,mcom,quo}. The
usual exchange statistics of quasiparticles is described a linear $B$ satisfying the
standard braid relations \begin{equation} \begin{array}{c} B_{(1)} B_{(2)} B_{(1)} =
B_{(2)} B_{(1)} B_{(2)} , \label{bra} \end{array} \end{equation} where $B_{(1)} := B
\otimes id$ and $B_{(2)} := id \otimes B$. The exchange process determined by the
operator $B$ is a real process. Such exchange process is possible if the dimension of
the effective space is equal or great than two. Hence in this case we need two
operators $T$ and $B$ for the description of our system with generalized statistics.
These operators are not arbitrary. They must satisfy the following consistency
conditions
\begin{equation}
\begin{array}{c}
B^{(1)}T^{(2)}T^{(1)} = T^{(2)}T^{(1)} B^{(2)},\\
(id_{E \otimes E} + \tilde{T})(id_{E \otimes E} - B) = 0,
\label{cod}
\end{array}
\end{equation}
where the operator $\tilde{T} : E\otimes E\rightarrow E\otimes E$ is given by its
matrix elements
\begin{equation} (\tilde{T})^{ij}_{kl} = T^{ki}_{lj}. \end{equation} We need a solution of these
conditions for the construction of an example of corresponding category.

\section{\hspace{-4mm}.\hspace{2mm}Hermitian Wick algebras
and Fock space representation}

Let us assume that a category of states $\mathcal{C}$ is given. We consider here a
pair of unital and associative algebras $\mathcal{A}$ and $\mathcal{A}^{\ast}$ in
$\mathcal{C}$. We assume that they are conjugated. This means that there is an
anti-linear and involutive anti-isomorphism $(-)^{\ast} :
\mathcal{A}\rightarrow\mathcal{A}^{\ast}$ and we have the following relations
\begin{equation} b^{\ast}\cdot a^{\ast} = (a \cdot b))^{\ast}, \quad (a^{\ast})^{\ast} = a, \end{equation}
where $a, b \in \mathcal{A}$ and $a^{\ast}, b^{\ast}$ are their images under the
ant-isomorphisms $(-)^{\ast}$. Both algebras $\mathcal{A}$ and $\mathcal{A}^{\ast}$
are graded
\begin{equation} \begin{array}{cc} \mathcal{A} := \bigoplus\limits_{n} \ \mathcal{A}^n ,& \mathcal{A}^{\ast} :=
\bigoplus\limits_{n} \ \mathcal{A}^{\ast n} . \end{array} \end{equation} A linear
mapping $\Psi :
\mathcal{A}^{\ast}\otimes\mathcal{A}\rightarrow\mathcal{A}\otimes\mathcal{A}^{\ast}$
such that we have the following relations
\begin{equation} \begin{array}{l} \Psi \circ (id_{\mathcal{A}^{\ast}} \otimes m_{\mathcal{A}}) = (m_{\mathcal{A}} \otimes
id_{\mathcal{A}^{\ast}})
\circ (id_{\mathcal{A}} \otimes \Psi) \circ (\Psi \otimes id_{\mathcal{A}}),\\
\Psi \circ (m_{\mathcal{A}^{\ast}} \otimes id_{\mathcal{A}}) = (id_{\mathcal{A}}
\otimes m_{\mathcal{A}^{\ast}}) \circ (\Psi \otimes id_{\mathcal{A}^{\ast}})
\circ (id_{\mathcal{A}^{\ast}} \otimes \Psi)\\
(\Psi(b^{\ast}\otimes a))^{\ast} = \Psi (a^{\ast} \otimes b) \label{twc} \end{array}
\end{equation} is said to be a cross symmetry or $\ast$--twist \cite{bma}. We use
here the notation
\begin{equation} \Psi(b^{\ast}\otimes a) = \Sigma a_{(1)}\otimes b^{\ast}_{(2)} \end{equation} for
$a\in\mathcal{A}, b^{\ast}\in\mathcal{A}^{\ast}$.

The tensor product $\mathcal{A}\otimes\mathcal{A}^{\ast}$ of algebras $\mathcal{A}$
and $\mathcal{A}^{\ast}$ equipped with the multiplication \begin{equation}
\begin{array}{c} m_{\Psi} := (m_{\mathcal{A}} \otimes m_{\mathcal{A}^{\ast}}) \circ (id_{\mathcal{A}} \otimes
\Psi \otimes id_{\mathcal{A}^{\ast}}) \label{mul}
\end{array}
\end{equation} is an associative algebra called a Hermitian Wick algebra
\cite{jswe,bma} and it is denoted by $\mathcal{W} = \mathcal{W}_{\Psi}(\mathcal{A}) =
\mathcal{A} \otimes_{\Psi} \mathcal{A}^{\ast}$. This means that the Hermitian Wick
algebra $\mathcal{W}$ is the tensor cross product of algebras $\mathcal{A}$ and
$\mathcal{A}^{\ast}$ with respect to the cross symmetry $\Psi$ \cite{bma}. Let $H$ be
a linear space. We denote by $L(H)$ the algebra of linear operators acting on $H$.
One can prove \cite{bma} that we have the

{\bf Theorem:} Let $\mathcal{W} \equiv \mathcal{A} \otimes_{\Psi}
\mathcal{A}^{\ast}$  be a Hermitian Wick algebra. If $\pi_{\mathcal{A}} : \mathcal{A}
\rightarrow L(H)$ is a representation of the algebra $\mathcal{A}$, such that we have
the relation
\begin{equation}
\begin{array}{c} (\pi_{\mathcal{A}}(b))^{\ast} \pi_{\mathcal{A}}(a) = \Sigma \pi_{\mathcal{A}}(a_{(1)})
(\pi_{\mathcal{A}}(b_{(2)}))^{\ast}, \label{wre}
\end{array}
\end{equation} then there is a representation $\pi_{\mathcal{W}} : \mathcal{W} \rightarrow L(H)$ of the
algebra $\mathcal{W}$. We use the following notation \begin{equation}
\pi_{\mathcal{A}}(x^i) \equiv a_{x^i}^+ , \quad \pi_{\mathcal{A}^{\ast}}(x^{\ast i})
\equiv a_{x^{\ast i}} .
\end{equation} The relations (\ref{wre}) are said to be {\it commutation relations}
if there is a positive definite scalar product on $H$ such that operators $a_{x^i}^+$
are adjoint to $a_{x^{\ast i}}$ and vice versa. Let us consider a Hermitian Wick
algebra $\mathcal{W}$ corresponding for a system with generalized statistics. For the
construction of such algebra we need a pair of algebras $\mathcal{A}$,
$\mathcal{A}^{\ast}$ and a cross symmetry $\Psi$. It is natural to assume that these
algebras have $E$ and $E^{\ast}$ as generating spaces, respectively, and there is the
following condition for the cross symmetry
\begin{equation} \Psi|_{E^{\ast}\otimes E} = T + g_E . \end{equation}

Let us consider the Fock space representation of the algebra $\mathcal{W}$
corresponding for a system with generalized statistics. For the ground state and
annihilation operators we assume that \begin{equation} \langle 0|0 \rangle = 0, \quad
a_{s^{\ast}} |0\rangle = 0 \quad \mbox{for} \quad s^{\ast} \in \mathcal{A}^{\ast}.
\end{equation} In this case the representation act on the algebra $\mathcal{A}$.
Creation operators are defined as the multiplication in the algebra $\mathcal{A}$
\begin{equation} a^+_{s} t := m_{\mathcal{A}}(s \otimes t), \quad \mbox{for} \quad s,
t \in \mathcal{A} .
\end{equation} The proper definition of the action of annihilation operators on the
whole algebra $\mathcal{A}$ is a problem.

If the action of annihilation operators are given in such a way that there is unique,
nondegenerate, positive definite scalar product on $\mathcal{A}$, creation operators
are adjoint to annihilation ones and vice versa, then we say that we have the
well--defined system with generalized statistics in the Fock representation
\cite{qstat}.

Let us consider some examples for such systems. Assume that quasiparticles and
quasiholes are moving on one dimensional effective space. In this case we can
construct the category $\mathcal{C} (E, E^{\ast}, g_{E})$ which consists of all
possible tensor product of spaces $E$ and $E^{\ast}$ and their direct sums.

The algebra of states $\mathcal{A}$ is the full tensor algebra $T E$ over the space
$E$, and the conjugate algebra $\mathcal{A}^{\ast}$ is identical with the tensor
algebra $TE^{\ast}$. If $T \equiv 0$ then  we obtain the most simple example of
well--defined system with generalized statistics. The corresponding statistics is the
so--called infinite Bolzman) statistics \cite{owg,gre,qstat}. The action of
annihilation operators is given by the formula \begin{equation} a_{x^{\ast i_k}
\otimes \cdots \otimes x^{\ast i_1}} (x^{j_1} \otimes \cdots \otimes x^{j_n}) :=
\delta_{i_1}^{j_1} \cdots \delta_{i_k}^{j_k} \ x^{j_{n-k+1}} \otimes \cdots \otimes
x^{j_n}.
\end{equation} For the scalar product we have the equation
\begin{equation} \langle i_n \cdots i_1 |j_1\cdots j_n \rangle^n := \delta^{i_1 j_1}\cdots \delta^{i_n
j_n}, \end{equation} where $i_n \cdots i_1 := x^{i_n}\cdots x^{i_1}$. It is easy to
see that we have the relation and \begin{equation} \begin{array}{c} a_{x^{\ast i}}
a^+_{x_j} := \delta_i^j {\bf 1}. \label{grel} \end{array} \end{equation}

Now let us consider the category $\mathcal{C} (E, E^{\ast}, g_{E}, |Psi^T)$ which
contains in addition a cross symmetry $\Psi^T$. If a cross operator $T :
E^{\ast}\otimes E\rightarrow E \otimes E^{\ast}$ is given, then the corresponding
cross symmetry $\Psi^T : TE^{\ast}\otimes TE\rightarrow TE\otimes TE^{\ast}$. is
defined by a set of mappings $\Psi_{k,l}:E^{\ast\otimes k}\otimes E^{\ast\otimes k}
\rightarrow E^{\otimes l}\otimes E^{\ast\otimes k}$, where $\Psi_{1,1} \equiv R := T
+ g_E$, and
\begin{equation}
\begin{array}{l}
\Psi_{1,l} := R^{(l)}_l \circ \ldots \circ R^{(1)}_l,\\
\Psi_{k,l} := (\Psi_{1,l})^{(1)}  \circ \ldots \circ (\Psi_{1,l})^{(k)}, \label{up}
\end{array} \end{equation} here $R^{(i)}_l : E^{(i)}_l \rightarrow E^{(i+1)}_l$, $E^{(i)}_l := E \otimes
\ldots \otimes E^{\ast} \otimes E \otimes \ldots \otimes E$ ($l+1$-factors,
$E^{\ast}$ on the i-th place, $i \leq l$) is given by the relation
$$
R^{(i)}_l := \underbrace{ id_{E} \otimes \ldots \otimes R \otimes \ldots \otimes
id_{E} }_{l\;\;times},
$$
where $R$ is on the i-th place, $(\Psi_{1,l})^{(i)}$ is defined in similar way like
$R^{(i)}$. The commutation relations (\ref{wre}) can be given here in the following
form \begin{equation} \begin{array}{c} a_{x^{\ast i}} a_{x^j}^+ - T^{ij}_{kl} \
a_{x^l}^+ a_{x^{\ast k}} = \delta^{ij}{\bf 1}. \label{crel} \end{array}
\end{equation} If the operator $\tilde{T}$ is a bounded operator acting on some
Hilbert space such that we have the following Yang-Baxter equation on $E\otimes E
\otimes E$
\begin{equation} (\tilde{T} \otimes id_E )\circ (id_E \otimes \tilde{T} )\circ (\tilde{T} \otimes id_E ) =
(id_E \otimes \tilde{T} ) \circ (\tilde{T}, \otimes id_E )\circ (id_E \otimes
\tilde{T} ),
\end{equation} and $||\tilde{T}|| \leq 1$, then according to Bo$\dot{z}$ejko and
Speicher \cite{bs2} there is a positive definite scalar product
\begin{equation}
\langle s | t \rangle^{n}_T := \langle s | P_n t \rangle^{n}_0 , \label{csca}
\end{equation}
where $s, t\in T E$,
$P_1 \equiv id$, $P_2 \equiv R_2 \equiv id_{E\otimes E} + \tilde{T}$, and
\begin{equation}
P_{n+1} := (id \otimes P_n) \circ R_{n+1},
\end{equation}
the operator $R_n$ is
given by the formula
\begin{equation} R_n := id + \tilde{T}^{(1)} + \tilde{T}^{(1)}
\tilde{T}^{(2)} + \cdots + \tilde{T}^{(1)}\dots\tilde{T}^{(n-1)} . \end{equation}
Note that the existence of nontrivial kernel of operator $P_2 \equiv id_{E\otimes E} +
\tilde{T}$ is essential for the nondegeneracy of the scalar product \cite{jswe}. One
can see that if this kernel is trivial, then we obtain well--defined system with
generalized statistics \cite{RM,ral}.

If the dimension of the effective space is great than one and the kernel of $P_2$  is
nontrivial, then the scalar product is degenerate. Hence we must remove this
degeneracy by factoring the mentioned above scalar product by the kernel. In this case
our category contains some quotients and we have $\mathcal{A} := TE/I, \quad
\mathcal{A}^{\ast} := TE^{\ast}/I^{\ast}$, where $I := gen\{id_{E\otimes E} - B\}$ is
an ideal in $T E$ and $B: E \otimes E \rightarrow E \otimes E$ is a linear and
invertible operator satisfying the braid relation (\ref{bra}) and the consistency
conditions (\ref{cod}), $I^{\ast}$ is the corresponding conjugated ideal in
$TE^{\ast}$. One can see that there is the cross symmetry and the action of
annihilation operators can be defined in such a way that we obtain the well--defined
system with the usual braid statistics \cite{RM,ral,quo}. We have here the following
commutation relations \begin{equation}
\begin{array}{l} a_{x^{\ast i}} a_{x^j}^+ - T^{ij}_{kl} \
a_{x^l}^+ a_{x^{\ast k}} = \delta^{ij}{\bf 1}\\
a_{x^{\ast i}} a_{x^{\ast j}} - B_{ij}^{kl} \
a_{x^{\ast l}} a_{x^{\ast k}} = 0, \\
a^+_{x^{i}} a^+_{x^j} - B^{ij}_{kl} \ a^+_{x^l} a^+_{x^{\ast k}} = 0. \label{brel}
\end{array}
\end{equation} Observe that for $T \equiv B \equiv \tau$, where $\tau$ represents the
transposition $\tau (x^i \otimes x^j) := x^j \otimes x^i$ we obtain the usual
canonical commutation relations of bosons or fermions. Note that similar relations are
described by Fiore \cite{fio}.

We can see that the commutation relations (\ref{crel}) corresponding for the system
with cross symmetry can be understood as a deformation of relations (\ref{grel}) of
Bolzman statistics. The commutation relations (\ref{brel}) for the system with a
braid group statistics are in fact certain degenerated system. In this way the
infinite statistics seems to be the most general statistics and others statistics are
its deformed or degenerated version.

\section{\hspace{-4mm}.\hspace{2mm}Category of comodules}

Let us consider a category of $H$--comodules, where $H$ is a finite Hopf algebra $H =
H(m, u, \triangle , \eta, S)$, equipped with the multiplication $m$, the unit $u$, the
comultiplication $\triangle$, the counit $\eta$ and the antipode $S$. According to our
physical interpretation the algebra represents an external quantum field and
comodules -- a charged particle dressed with quanta of this field. We use the
following notation for the coproduct in $H$: if $h \in H$, then $\triangle (h) :=
\Sigma h_{_{(1)}} \otimes h_{_{(2)}} \in H \otimes H$. We assume that $H$ is
coquasitriangular Hopf algebra (CQTHA). This means that $H$ equipped with a bilinear
form $b : H \otimes H \rightarrow \mathbb{C}$ such that \begin{equation}
\begin{array}{l} \Sigma b(h_{_{(1)}}, k_{_{(2)}}) k_{_{(2)}} h_{_{(2)}}
= \Sigma h_{_{(1)}} k_{_{(1)}} b(h_{_{(2)}}, k_{_{(2)}}),\\
b(h, kl) = \Sigma b(h_{_{(1)}}, k) b(h_{_{(2)}}, l),\\
b(hk, l) = \Sigma b(h, l_{_{(2)}}) b(k, l_{_{(1)}}) \end{array} \end{equation} for
every $h, k, l \in H$. If such bilinear form $b$ exists for a given Hopf algebra $H$,
then we say that there is a {\it coquasitriangular structure} on $H$. Let us assume
for simplicity that $H \equiv \mathbb{C} G$ is a group algebra, where $G$ is an
Abelian group. The group algebra $H := \mathbb{C} G$ is a Hopf algebra for which the
comultiplication, the counit, and the antipode are given by the formulae
$$
\begin{array}{cccc} \triangle (g) := g \otimes g,&\eta(g) := 1,&S(g) := g^{-1}&\mbox{for} \ g \in G. \end{array}
$$
respectively. If $G$ is an Abelian group, then the coquasitriangular structure on $H
= \mathbb{C} G$ is given by a commutation factor $\epsilon : G\otimes
G\rightarrow\mathbb{C}\setminus\{0\}$ on $G$, \cite{mon}.

A right $H$-Hopf module is a $k$-linear space $M$ such that\\
(i) there is a right $H$-module action $\lhd : M \otimes H \rightarrow M$,\\
(ii) there is a right $H$-comodule map $\delta : M \rightarrow M \otimes H$, \\
(iii) $\delta$ is a right $H$-module map, this means that we have the relation
\begin{equation} \Sigma \ (m \lhd h)_{(0)} \otimes (m \lhd h)_{(1)} = \Sigma \ m_{(0)} \lhd
h_{(1)} \otimes m_{(1)} \ h_{(2)}, \end{equation} where $m \in M, h \in H, \delta (m)
= \Sigma \ m_{(0)} \otimes m_{(1)}$, and $\triangle (h) = \Sigma \ h_{(1)} \otimes
h_{(2)}$.

It is very interesting that every right $H$-Hopf module $M$ is a tensor product $E =
\mathcal{V} \otimes H$, where $\mathcal{V}\equiv E^{coH}$ is also a trivial right
$H$-module, i. e. $m \lhd h = \eta (h) m$ for every $m \in E$. If $\mathcal{V}$ is a
finite--dimensional vector space equipped with a basis $\{\xi^i : i=1,\ldots,n\}$ and
$H$ is a finite Hopf algebra equipped with a vector space basis $\{h^i :
i=1,\ldots,N\}$, then we have the following basis in $E$ \begin{equation} x^{i} :=
\xi^i \otimes h^i \quad\mbox{(no sum)}.
\end{equation} In this basis the right coaction is given by the relation
\begin{equation} \begin{array}{c} \delta (x^i) := \Sigma \ x^{i (1) } \otimes h^{i (2)} . \label{hco}
\end{array} \end{equation}

Let $H$ be a CQTHA with coquasitriangular structure $( -,-)$. The family of all
$H$-comodules forms a category $\mathcal{C} = \mathcal{M}^H$. The category
$\mathcal{C}$ is braided monoidal. The braid symmetry $\Psi \equiv \{\Psi_{U, V} :
U\otimes V\rightarrow V\otimes U; U, V\in Ob\mathcal{C}\}$ in $\mathcal{C}$ is
defined by the equation
\begin{equation}
\begin{array}{c} \Psi_{U, V} (u \otimes v) = \Sigma b(v_1 , u_1) \ v_0 \otimes u_0 ,
\label{coin}
\end{array}
\end{equation} where $\rho (u) = \Sigma u_0 \otimes u_1 \in U \otimes H$, and $\rho (v) =
\Sigma v_0 \otimes v_1 \in V \otimes H$ for every $u \in U , v \in V$. Let $H$ be a
CQTHA with coquasitriangular structure given by a form $b : H\otimes H \rightarrow k$
and let $\mathcal{A}$ be a (right) $H$-comodule algebra with coaction $\rho$. Then the
algebra $\mathcal{A}$ is said to be {\it quantum commutative} with respect to the
coaction of $(H, b)$ if an only if we have the relation
\begin{equation} \begin{array}{c} a \ b = \Sigma \ b(a_1, b_1) \ b_0 \ a_0 , \end{array} \end{equation}
where $\rho (a) = \Sigma a_0 \otimes a_1 \in \mathcal{A} \otimes H$, and $\rho (b) =
\Sigma b_0 \otimes b_1 \in \mathcal{A} \otimes H$ for every $a, b \in \mathcal{A}$,
see \cite{cowe}. The Hopf algebra $H$ is said to be a {\it quantum symmetry} for
$\mathcal{A}$.

If $H \equiv kG$, where $G$ is an Abelian group, $k \equiv \mathbb{C}$ is the field of
complex numbers, then the coquasitriangular structure on $H$ is given as a
bicharacter on $G$ \cite{mon}, i. e. a mapping $\epsilon : G \times G \rightarrow
\mathbb{C} \setminus \{0\}$ such that we have the following relations
\begin{equation} \epsilon (\alpha, \beta + \gamma) = \epsilon (\alpha, \beta)
\epsilon (\alpha , \gamma), \quad \epsilon (\alpha + \beta , \gamma) = \epsilon
(\alpha , \gamma) \epsilon (\beta , \gamma)
\end{equation} for $\alpha, \beta, \gamma \in G$. If in addition \begin{equation} \epsilon (\alpha, \beta) \epsilon (\beta, \alpha)
= 1, \end{equation} for $\alpha, \beta \in G$. It is interesting that the coaction of
$H$ on certain space $\mathcal{V}$, where $H \equiv k G$ is equivalent to the
$G$-gradation of $\mathcal{V}$,see \cite{mon}. In this case the quantum commutative
algebra $\mathcal{A}$ becomes graded commutative \begin{equation} x^i x^j = \epsilon
(\mid i\mid , \mid j\mid) x^j x^i
\end{equation} for homogeneous elements $x^i$ and $x^j$ of grade $\mid i\mid$ and $\mid j\mid$,
respectively. One can describe the corresponding algebra $\mathcal{A}^{\ast}$ in a
similar way. For both two operators $T$ and $B$ we have the following relations
\begin{equation} T(x^{\ast i}\otimes x^j) := \epsilon (\mid j\mid , \mid i\mid) x^j\otimes x^{\ast
i}, \quad B(x^i\otimes x^j) := \epsilon (\mid i\mid , \mid j\mid) x^j\otimes x^i
\end{equation} The construction of related Hermitian Wick algebra and Fock space
representation is evident, see \cite{quo}. Note that in this case the category of
comodules $\mathcal{C} = \mathcal{M}^H$ is denoted by $\mathcal{C} \equiv
\mathcal{C}(G, \epsilon )$.

\section*{Acknowledgments}
The author would like to thank A. Borowiec for the discussion and any other help.

\end{document}